\documentclass[12pt]{amsart}
\usepackage{amsmath,amssymb, euscript, graphicx, hyperref}

\newtheorem{thm}{Theorem}[section]
\newtheorem{lem}[thm]{Lemma}
\newtheorem{prop}[thm]{Proposition}
\theoremstyle{definition}
\newtheorem{defn}[thm]{Definition}

\newtheorem{rem}[thm]{Remark}
\newtheorem{cor}[thm]{Corollary}

\newcommand{\blackboard}[1]{\ensuremath{\mathbb{#1}}}

\newcommand{\N}{\blackboard{N}}

\newcommand{\Z}{\blackboard{Z}}

\newcommand{\Q}{\blackboard{Q}}
\newcommand{\R}{\blackboard{R}}
\newcommand{\C}{\blackboard{C}}

\address{Azer Akhmedov, Department of Mathematics,
North Dakota State University,
Fargo, ND, 58102, USA}
\email{azer.akhmedov@ndsu.edu}

 \address{Damiano Fulghesu, Department of Mathematics,
Minnesota State University,
Moorhead, MN, 56563, USA}
\email{damiano.fulghesu@mnstate.edu}

\begin{document}

\begin{center} {\bf \LARGE Arithmetic Sets in Groups} \end{center}

\medskip

\begin{center} AZER AKHMEDOV and DAMIANO FULGHESU \end{center}

\bigskip

{\Small Abstract: We define a notion of an arithmetic set in an arbitrary countable group and study properties of these sets in the cases of Abelian groups and non-abelian free groups.} 

\section{Introduction}

 To motivate the notion of an arithmetic set in the group $\mathbb{Z}$ we would like to consider the following two systems  of equations over the reals where $M$ is a given real number:
 
 $$x_{n}+x_{n+1}+x_{n+2} = M, n\in \mathbb{Z} \ \ (A) $$
 
 and 
 
 $$x_{n}+x_{n+1}+x_{n+3} = M, n\in \mathbb{Z} \ \ (B). $$  
 
 \medskip
 
 Both of the systems admit a non-degenerate (i.e. nonconstant) solution for any real $M$. Let us observe that any solution $${\bf x} = (\dots , x_{-1}, x_0, x_1, \dots )$$ of the system (A) satisfies the following properties:
 
 \medskip
 
 (A1) ${\bf x}$ is periodic, i.e. there exist $a, b, c\in \mathbb{R}$ such that $x_{3k} = a, x_{3k+1} = b, x_{3k+2} = c$ for all integers $k$;   
 
 \medskip
 
 (A2) ${\bf x}$ is bounded (since it is periodic).

 \medskip
 
 In sharp contrast, for the solutions of the system (B) the following properties hold:
 
 \medskip 
 
 (B1) any periodic solution is degenerate;
 
 \medskip
 
 (B2) any bounded solution is degenerate.
  
 \medskip
 
 It is immediate to see property (B1), while a stronger property (B2) follows from the fact that if ${\bf x} = (\dots , x_{-1}, x_0, x_1, \dots )$ is a solution then $x_n = c_1\lambda _1^n + c_2\lambda _2^n + c_3\lambda _3^n + c_4$ for any $n\geq 1$ where $\lambda _1, \lambda _2, \lambda _3$ are the three non-identity roots of the polynomial $p(z) = z^4-z^3+z^2-1$ and  $c_1, c_2, c_3, c_4$ are some complex numbers. Here, one needs to observe that the polynomial $p(z)$ has four distinct roots $\lambda _1, \lambda _2, \lambda _3, \lambda _4$ where $\lambda _4 = 1$, $\lambda _3\in (-1, 1)$ while the non-real complex conjugate roots $\lambda _1$ and $\lambda _2$ lie outside the unit circle in the complex plane.  
 
 \medskip
 
 One might ask a question: why is that the properties of the systems (A) and (B) are so drastically different? Let us also observe that the converse of the property (A1) also holds, namely, for any real $a, b, c$ with $a+b+c=M$ one immediately obtains a solution of the system (A) by letting $x_{3k} = a, x_{3k+1} = b, x_{3k+2} = c$ for all integers $k$. (On the other hand, no such simple algorithm exists for obtaining the solution of the system (B)). 
 
 \medskip
 
 Notice that both of the systems (A) and (B) can be written as $$\displaystyle \sum _{g+n\in S}x_g = M, \ n\in \mathbb{Z}(\star )$$ where $S$ is a finite subset of $\mathbb{Z}$; for the system (A) we will have $S = \{0,1,2\}$ while for the system (B) we take $S = \{0,1,3\}$. Then the striking difference of the properties of (A) and (B) can perhaps be explained by the fact that the set $\{0,1,2\}$ is a tile of $\mathbb{Z}$ while the set $\{0,1,3\}$ is not!

  \bigskip
  
  In this paper, we investigate the above asked questions, i.e. how the combinatorics of the finite subset $S$ influences the properties of the solution of the system $(\star )$. We find out that the issue is more subtle than $S$ just being (or not being) a tile. By our studies we are lead to the notion of an arithmetic subset (in $\mathbb{Z}$, and more generally, in an arbitrary countable group).
  
  \medskip
  
   While no complete characterization of the tiles of $\mathbb{Z}$ is available, it is well known that being a tile imposes strong conditions of arithmetic flavor on the structure of a set (e.g. see \cite{LW, N, K, RX, S}). The Coven-Meyerowitz Conjecture states that being a tile is equivalent to a purely arithmetic condition on the set (see below, in Section 2). It is interesting to examine the arithmetic condition proposed by Coven-Meyerowitz Conjecture with respect to the notion of arithmeticity of the set studied in this paper. At the end of Section 4, we point out some relevant properties of arithmeticity related to the Coven-Meyerowitz Conjecture.     
  
 \medskip
 
  In the current paper, for an arbitrary (countable) group, we introduce the notion of an arithmetic subset. Every tile of a group turns out to be arithmetic while arithmetic sets form a larger class of subsets. It turns out that, in negatively curved groups (such as free groups of rank at least two) being arithmetic is a more loose condition on a set while in groups at the other extreme (such as Abelian or nilpotent groups) it imposes very strong conditions of number-theoretic flavor. For more in depth study of the notion of arithmeticity, we limit ourselves primarily to two examples of groups: free groups and free Abelian groups. 
  
  \medskip
  
  {\em Structure of the paper:} In Section 2, we review basic facts about tiles of groups. In Section 3, we introduce the major notions of the paper and observe simple properties. In Sections 4 and 5, we study properties of arithmetic sets in Abelian groups and free groups respectively.
   
\bigskip
 
\section{Tiles in Groups}

\begin{defn} Let $\Gamma $ be a countable group. A finite set $K\subseteq \Gamma $ is called a {\em right tile} if there exists a subset $C\subseteq \Gamma $ such that $\Gamma = \displaystyle \mathop{\sqcup }_{g\in C}gK$. 
\end{defn}

 Similarly, one can introduce left tiles of groups as finite subsets whose right shifts form a tiling, however, in this paper we will be considering only the right tiles. Notice that a subset $K$ is a right tile if and only if $K^{-1}$ is a left tile. Obviously, for Abelian groups, the two notions coincide.
 
 \medskip

  {\bf Example 1.}  The set $\{x,y\}$ is a tile of a group $\Gamma $ if and only if the element $x^{-1}y$ is either non-torsion or has even order.

  \medskip

  {\bf Example 2.}  The set  $\{0,1,3\}$ \ is not a tile of $\mathbb{Z}$.  More generally, the set \ $\{0, 1, x\}$ is a tile if and only if $x\equiv 2 \ (\mathrm{mod}3)$. On the other hand, any finite arithmetic progression in $\mathbb{Z}$ is a tile of $\mathbb{Z}$.

  \medskip

  {\bf Example 3.} Let $p$ be a prime number, and $$K = \{a_0 , a_1, a_2, \dots , a_{p-1}\}\subset \mathbb{Z}$$ be a finite subset where $a_0 < a_1 < \dots < a_{p-1}$. Then $K$ is a tile if and only if there exists $k\in \mathbb{N}$ such that the cyclotomic polynomial $\Phi _{p^k}(z)$ divides $P_K(z) = \displaystyle \sum _{i=1}^{p}z^{b_i}$ where $b_i = a_i - \mathrm{min} K$ for all $i\in \{0, 1, \ldots , p-1\}$. This is a rephrasement of the result from \cite{N} which states that $K$ is a tile if and only if for some non-negative $e$, all elements of $K$ are congruent mod$(p^e)$ and incongruent mod$(p^{e+1})$. 
  
  \medskip
  
  Newman (\cite{N}) has determined all tiles $K$ of $\mathbb{Z}$ where $|K|$ is a prime power. In general, however, the problem remains open. A positive solution to the following conjecture would yield a nice characterization of tiles in $\mathbb{Z}$.

  \medskip

  {\bf Conjecture of Coven-Meyerowitz, [\cite{CM}, 1999] :}   Let $K\subset \mathbb{Z}$ be a finite subset, $R_K = \{ d\in \mathbb{N} \ | \ \Phi _d \ divides \ P_K(x)\}$, and $S_K = \{ p^{\alpha }\in R_K\}$ - the set of prime powers of $R_K$. Then $K$ is a $\mathbb{Z}$-tile if and only if the following conditions are satisfied:
  
  \medskip

   $(T_1)$ : $P_K(1) = \displaystyle \mathop{\Pi }_{s\in S_K}\Phi _s(1)$

   $(T_2)$ :  if $x_1, x_2, \ldots , x_n \in S_K$ are powers of distinct primes then $x_1x_2\ldots x_n \in R_K$.
   
   \medskip

   It is known that $(T_1) \ and \ (T_2) \Rightarrow K \ is \ a \ \mathbb{Z}$-$tile$. Moreover, $K \ is \ a \
   \mathbb{Z}$-$tile \Rightarrow (T_1)$. It is not known whether or not $K \ is \ a \ \mathbb{Z}$-$tile \Rightarrow (T_2)$

   \medskip

  {\bf Example 4.} If $\mathbb{Z}^2 = \langle a, b \rangle $ then the set $\{a^ib^j | -n \leq i, j \leq n\}$ is an obvious  tile for any $n\in \mathbb{N}$. This set can be viewed as a ball of $\mathbb{Z}^2$ with respect to the generating set $S = \{a^{\pm 1}, b^{\pm 1}, a^{\pm 1}b^{\pm 1}, a^{\pm 1}b^{\mp 1}\}$. On the other hand, it is easy to find finite symmetric generating sets of $\mathbb{Z}^2$ with respect to which the balls of positive radii are not tiles.

  \medskip

  For general $\mathbb{Z}^n$, we would like to mention the following

  \medskip

  {\bf Fuglede Conjecture} [\cite{F}, 1974] : If $A\subseteq \mathbb{R}^n$ is a measurable subset then $A$ is a tile of $\mathbb{R}^n$ if and only if $A$ is spectral, i.e. for some set (spectrum) $\Lambda \subseteq \mathbb{R}^n$, \ the space $L^2(A)$ has an orthogonal basis $\{e^{2\pi i\lambda x}\}_{\lambda \in \Lambda }$.

\medskip

  The conjecture has been disproved by T.Tao \cite{T} for $n\geq 5$, and \cite{M}, \cite{FMM} extended the result to the cases of  $n = 4$ and $n=3$ respectively. For $n = 1, 2$ it still remains open.

  \medskip

  {\bf Example 5.} Let $r\in \mathbb{N}$ and $B_r$ be a ball of radius $r$ in the Cayley graph of the free group $\mathbb {F}_k = \langle a_1, \ldots , a_k \rangle $ with respect to standard generating set. Then the sets $B_1$ and $B_1\backslash \{1\}$ are tiles of $\mathbb {F}_k$. On the other hand, it is not difficult to see that $B_r$ is a tile for every $r\geq 2$ as well (see Proposition \ref{prop:prop1}), while $B_r\backslash \{1\}$ is not a tile for $r\geq 2$.

  \medskip

  {\bf Example 6.} If $\mathbb {F}_k = \langle a_1, \ldots , a_k \rangle $ and $k\geq 2$ then a sphere $S_r$ in the Cayley graph w.r.t. standard generating set is not a tile if $r\geq 1$ is even.

  \bigskip

\section{Arithmetic Sets}

Let $G$ be a finitely generated group and let $K=\{ g_0, \dots, g_{k-1}\} \subset G$ be a proper subset of $G$ which generates $G$. We associate to $K$ the system of equations:
$$
A(K)= \left\{ \sum_{i=0}^{k-1} x_{g_i \cdot g} = 0 \right\}_{g \in G}
$$
In case $G$ is finite, $A(K)$ will also represent the corresponding square $(0,1)$-matrix. Clearly, for every $G$ and $K$, we have the trivial solution $\{ x_g = 0\}_{g \in G}$.

\medskip

\begin{defn}
A solution $\{ x_g = \alpha_g\}_{g \in G}$ for the system $A(K)$ is called {\em bounded} if there exists a real number $M$ such that, for all $g \in G$, we have $\left| \alpha_g \right| \leq M$. 

\medskip

A solution $\{ x_g = \alpha_g\}_{g \in G}$ for the system $A(K)$ is called {\em periodic} if there exists a finite subset $J$ and a subset $I$ of $G$ such that the following two conditions hold:

 (i) $\displaystyle \bigcup_{h \in I} \{ h \cdot s \}_{s \in J} = G$ and the intersection of every pair of sets in the union is empty, in particular the finite set $J$ is a right tile. More precisely, the lateral classes $\{ Is\}_{s \in J}$ are a partition of $G$. 
 
  (ii) $\alpha_{h_1 \cdot s}=\alpha_{h_2 \cdot s}$ for every pair of elements $h_1$ and $h_2$ in $I$.
\end{defn}

\medskip

\begin{defn}
Let $G$ and $K$ be as above. We say that $K$ is {\em b-arithmetic} in $G$ if the system $A(K)$ has a non-trivial bounded solution in $\C$. We say that $K$ is {\em p-arithmetic} in $G$, if the system $A(K)$ has a non-trivial periodic solution in $\C$. We say that $K$ is {\em totally b-arithmetic} in $G$ (respectively, {\em totally p-arithmetic}) if all solutions are bounded (respectively periodic).
\end{defn}

\medskip

\begin{rem}
Clearly, if $K$ is p-arithmetic in $G$, then it is also b-arithmetic in $G$. Moreover, if $G$ is finite any solution is bounded and periodic, therefore, we add the following definition.
\end{rem}

\medskip

\begin{defn}
If $G$ is finite, we say that $K$ is arithmetic in $G$ if and only if $A(K)$ has a non trivial solution.
\end{defn}

\medskip

\begin{thm}\label{tile.implies.arithmetic}
If a proper finite subset $K \subset G$ is a right tile in $G$, then $K$ is p-arithmetic in $G$.
\end{thm}

\medskip

{\bf Proof.} Since $K$ is a right tile, there exists a subset $I \subset G$ such that:
$$
\bigcup_{h \in I} \{ h \cdot g_i \}_{i=0 \dots k-1} = G
$$
and the intersection of every pair of sets in the above union is empty. We write, for every $g \in G$, $x_{g} = k-1$ if $g^{-1} \in I$ and $x_g=-1$ otherwise. We want to prove that $\{ x_{g}\}_{i=0 \dots k-1}$ is a solution of $A(K)$ (it is clear that such solution is periodic). It is enough to show that, for every $g \in G$, the set $\{ g_i \cdot g\}_{i=0 \dots k-1}$ contains exactly one element $g_i \cdot g$ such that $(g_i \cdot g )^{-1} \in I$. Now, for all $g$ in $G$, there exists a unique $g_i$ in $K$ such that $g^{-1} \in I g_i$. This happens if and only if $(g_i \cdot g )^{-1} \in I$.

\bigskip

\section{Arithmetic Subsets of $\Z$}

Without loss of generality, we will assume 
$$
K=\{ s_0, s_1, s_2, \dots, s_{k-1}\}
$$
where the $s_i$ are integers such that $0 =s_0< s_1 < s_2 < \dots < s_{k-1}$. Moreover $K$ will generate $G=\Z$, that is to say $s_1, \dots, s_{k-1}$ are relatively prime.

\medskip

For any initial choice of $s_{k-1}$ numbers $x_0, x_1, \dots, x_{s_{k-1}-1} $, which can be chosen in $\C$, $\R$, or $\Q$, we define a sequence by recurrence:
$$
x_{n}=\left\{ \begin{array}{l} -\sum_{i=0}^{k-2}x_{n - s_{k-1} + s_i} \text{ if  $n\geq s_{k-1}$}\\ \; \\ -\sum_{i=1}^{k-1}x_{n + s_i} \text{ if  $n<0$}\end{array} \right.
$$

Clearly the sequence $\left \{ x_i \right \}_{i \in \Z}$ is a solution for the system $A(K)$.

\medskip

We define also the mask polynomial 
$$
P_K(x)=\sum_{i=0}^{k-1} x^{s_i}
$$
It is known from the theory of recurrence sequences, that, if $\alpha_1, \dots, \alpha_{s_{k-1}}$ are roots of $P_K(x)$, then the generic term of the sequence can be written as
$$
x_n=\sum_{j=1}^{s_{k-1}}b_j \alpha_j^n
$$
where the coefficients $b_j$ depend on the first $s_{k-1}$ terms of the sequence. 

\bigskip

\begin{prop} \label{thm:bounded}
The sequence $\{ x_n\}_{n \in \Z}$ is bounded if and only if the coefficients $b_j$ corresponding to roots whose modulus is different from 1, are 0.
\end{prop}

{\bf Proof.} If all the roots have modulus one then there is nothing to prove.

\medskip

Up to reordering the roots, assume that $\{\alpha_1, \dots, \alpha_t\}$ is the set of roots of $P_K(x)$ with modulus greater than one. Then $x_n$ is bounded if and only if the sequence  $$y_n:=b_1\alpha^n_1 + \dots + b_t \alpha^n_t $$ is bounded as $n$ varies in $\N$. We can also assume that all coefficients $b_j$ are different from 0. Without loss of generality, let us assume that the largest modulus is attained by $\alpha_1, \dots , \alpha _r$ where $r\leq t$. Then for the proof of boundedness of $(y_n)_{n\geq 1}$ it suffices to show that the sequence $$z_n:=b_1\alpha^n_1 + \dots + b_r \alpha^n_r$$
is bounded as $n$ varies in $\N$. Then we can write
$$ z_n =\rho ^n(b_1\theta^n_1 + \dots + b_r\theta^n_r ) $$

 where $\rho = |\alpha _1| = \dots = |\alpha _r|$. Now, by compactness of the unit torus $\mathbb{T}^r$, the orbit $(\theta _1^n, \dots , \theta _r^n)_{n\geq 1}$ has an accumulation point $(\omega _1, \dots , \omega _r)$. 
 Then, since $\rho > 1$, the sequence $z_n$ is bounded only if $b_1\omega _1^{n} \dots + b_r\omega _r^{n} = 0$ for infinitely many $n$. Then, using the Wandermont determinant, one obtains that $b_1 = \dots = b_r = 0$. 
 
 \medskip
  
 Thus we proved that the only way to have the sequence $y_n$ bounded is that some of the coefficients $b_1, \dots , b_t$ are zero, getting a contradiction. 
  
  \medskip
  
 On the other hand, if all roots have moduli less or equal to one, we consider the sequence

$$
x_n=\sum_{j=1}^{s_{k-1}}b_j \alpha_j^{-n}
$$

and argue as above. 

\medskip

 Now it remains to consider the case when not all roots of the polynomial $P_K(x)$ are distinct. But in this case, instead of $y_n$ we need to consider the expression $y_n' = b_1n^{p_1}\alpha^n_1 + \dots + b_t n^{p_t}\alpha^n_t$ where $p_1, \dots , p_t$ are some non-negative integers. Then (instead of $z_n$) we will have $z_n' =  b_1n^{p_1}\alpha^n_1 + \dots + b_r n^{p_r}\alpha^n_r$. Here, let $p = \max \{p_1, \dots , p_r\}$. Again, up to reordering the indices, we obtain that the sequence $z_n'$ is bounded if and only if the sequence $z_n'' = b_1n^{p}\alpha^n_1 + \dots + b_s n^{p}\alpha^n_s$ is bounded where $s\leq r$. Then we obtain a contradiction as in the previous case. $\square $

\medskip

  In addition, let us point out that if the mask polynomial $P_K(x)$ vanishes for a complex number $\alpha$ whose modulus is one, since its coefficients are rational, $P_K(x)$ is divisible by the minimal polynomial for $\alpha$, whose roots have all modulus 1. Therefore, the following definition makes sense.

\begin{defn}
Let $K$ be as above, and let $P_K(x)$ be its mask polynomial. We will consider the decomposition
$$
P_K(x)=C_K(x) \cdot D_K(x)
$$
where $C_K(x)$ is either 1 or a product of polynomials whose roots have modulus one and $D_K(x)$ is a polynomial whose roots have modulus different from 1.
\end{defn}

The following Proposition is now clear.

\begin{prop}\label{thm:modulus}
A finite subset $K \subset \Z$ is b-arithmetic in $\Z$ if and only if $$C_K(x) \neq 1.$$
\end{prop}

We now focus on arithmetic sets in $\Z$ which have a periodic solution. In particular we have the following Proposition.

\begin{thm}\label{thm:equivalent}
Let $K$ be a subset of $\Z$ as above. Then the following are equivalent.

  a) $A(K)$ has a periodic solution;
  
  b) $K$ is arithmetic in $\Z_n$ for some $n>k=|K|$;
 
  c)  $P_K(x)$ is divisible by a cyclotomic polynomial $\Phi_m(x)$, such that $m|n$ where $n$ is as part (b); 

  d)  $A(K)$ has an integral periodic solution.

\end{thm}

{\bf Proof.}

 (a) implies (b). Let $p$ be the period of a periodic solution $\{ x_i \}_{i \in \Z}$. Fix $n:=lcm(p,k)$, then $\{ x_i \}_{i=0, \dots, n-1}$ is a solution for $A(K)$ in $\Z_n$.

 (b) implies (c). First of all notice that the matrix $A$ associated to $K$ in $\Z_n$ is circulant. From \cite{MM} Section I.4.9, we know that the eigenvalues of $A$ are exactly
$$
P_K(1), P_K(\xi_n), P_K(\xi_n^2), \dots, P_K(\xi_n^{n-1})
$$
where $\xi_n=e^{2\pi i/n}$ is a primitive $n$th-root of unit. Therefore $A$ is singular if and only if $P_K(\xi_n^{d})=0$ for some $d<n$. Now, let $h=\text{GCD}(d,n)$, we have that $\xi_n^{d}$ is a primitive root of order $m=n/h$. Since all the coefficients of $P_K(x)$ are 1 (hence rational), we have that $P_K(x)$ is divisible by $\Phi_m(x)$. By construction $m|n$.

 (c) implies (b). Assume that $P_K(x)$ is divisible by a cyclotomic polynomial $\Phi_m(x)$ and fix $n$ such that $n$ is a multiple $m$ larger than $k$. Consider the matrix $A$ associated to $K$ in $\Z_n$. Let $\xi_n=e^{2\pi i/n}$, we have $P_K \left( \xi^{n/m}_n\right)=0$, therefore, again from \cite{MM} Section I.4.9, the determinant of $A$ is 0.

(b) implies (d). Since the determinant of the matrix $A$ associated to $K$ in $\Z_n$  is 0 and all its coefficients are in $\Q$, then there is a solution $\{ x_1, \dots, x_n \}$ in $\Q$ for $\Z_n$. Now, take the least common denominator $d$ of all $x_i$ and multiply the solution by $d$. We get an integral solution. We then extend it by periodicity on all $\Z$.

 (d) implies (a). This is obvious. $\square $

\begin{lem}\label{k3cyclotomic}
Assume that $k=|K|=3$ and $K$ generates $\Z$, then $K$ is b-arithmetic in $\Z$, if and only if $P_K(x)$ is divisible by $\Phi_3(x)$. 
\end{lem}

{\bf Proof.}
Because of our initial assumptions on $K$ we have
$$
P_K(x) = 1 + x^a + x^b
$$
for some distinct integers $a$ and $b$, moreover, since $K$ generates $\Z$, we also have that $a$ and $b$ are relatively prime. From Proposition \ref{thm:modulus}, $K$ is b-arithmetic in $\Z$, if and only if there exists a complex number $\xi$ such that $P_K(\xi)=0$ and $|\xi|=1$. In particular $v:=\xi^a$ and $w:=\xi^b$ are two complex numbers on the unit circle whose sum is $-1$. The fact that the sum of $v$ and $w$ is real, implies that their imaginary parts are opposite. This, together with  the fact that these two numbers are on the unit circle, implies that the absolute values of their real parts are equal. In conclusion we must have
$$
\{ v, w \} = \{ \omega, \omega^2 \}
$$
where $\omega$ is a primitive cubic root of 1, therefore $\xi^{3a}=\xi^{3b}=1$ and $\xi^a, \xi^b \neq 1$. This means that $\xi$ is a primitive $m$-th root of 1 such that $m|3a$ and $m|3b$. Since $a$ and $b$ are relatively prime, we must have $m=3$. Moreover we must have that $a$ and $b$ are, in some order, $1$ and $2$ modulo 3. On the other hand, if $a$ and $b$ are congruent to $1$ and $2$ modulo 3, then the polynomial
$$
1 + x^a + x^b =0
$$
is divisible by $\Phi_3(x)$.

From Lemma \ref{k3cyclotomic} and Theorem \ref{thm:equivalent} we get the following statement.

\begin{cor}\label{equivalencek3}
Assume $k=3$, then the following are equivalent

a)  $K$ is b-arithmetic in $\Z$;

b) $A(K)$ has an integral periodic solution;

c) $K$ is arithmetic in $\Z_n$ for some $n>k=|K|$.
\end{cor}

\begin{prop}\label{corolk3}
A proper subset $K$, such that $|K|=3$, is b-arithmetic in $\Z$, if and only if $K$ is a tile of $\Z$.
\end{prop}

{\bf Proof.} From the proof of Lemma \ref{k3cyclotomic}, we have that $K=\{ 0, a, b \}$ is b-arithmetic in $\Z$, if and only if $a$ and $b$ are congruent to $1$ and $2$ modulo 3. But, by the result of Newmann \cite{N}, this happens if and only if $K$ is a tile. $\square $

\begin{prop}
For every integer $k>3$, there exists a finite set $K \subset \Z$ of cardinality $k$, such that $K$ is p-arithmetic, but not a tile for $\Z$.
\end{prop}

{\bf Proof.} We split the proof into two cases: $k$ is a prime greater than 3, and $k$ is not prime.

Assume $k>3$ is a prime number $p$. Consider the polynomial
$$
P(x):=\Phi_6(x)\cdot \left( 1 + x + \sum_{i=0}^{p-3} x^{p+i} \right)=(1-x +x^2)\cdot \left( 1 + x + \sum_{i=0}^{p-3} x^{p+i} \right).
$$

 A straightforward computation shows that

$$
P(x) = 1 + x^3 + x^p + \sum_{i=2}^{p-3} x^{p+i} + x^{2p-1}.
$$

Now, the set $K:= \{0, 3, p, p+2, \dots, 2p -3, 2p-1\}$, is p-arithmetic in $\Z$ because its mask polynomial $P(x)$ is divisible by a cyclotomic polynomial $\Phi_6$ (see Theorem \ref{thm:equivalent}).

On the other hand, $K$ is a tile only if there exists some non negative integer $e$ such that all the elements of $K$ are congruent (mod $p^e$) and all incongruent (mod $p^{e+1}$). Now $e$ cannot be zero because the elements of $K$ are not all incongruent (mod $p$), similarly $e$ cannot be greater than 0 because the elements of $K$ are not all congruent (mod $p^e$)  therefore $K$ cannot be a tile.

Now, assume that $k$ is not prime. Let us write $k=pd$ for some prime $p$ and $d>1$. We define $K$ as the difference of two sets 
$$
K:=\{ 0, 1, \dots, dp + d - 1\}\backslash \{ p, 2p + 1, 3p + 2, \dots, dp + d - 1\}.
$$
Clearly $P_K(x)$ is divisible by $\Phi_p(x)$, therefore (again by Theorem \ref{thm:equivalent}) $K$ is p-arithmetic. On the other hand, $K$ cannot be a tile because any left shift of $K$ covering $p$ will necessarily overlap $K$. $\square $ 

\bigskip

\begin{prop}
Let $K$ be a finite set in $\Z$ such that $|K|=4$. Then $K$ is p-arithmetic if and only if $K$ is b-arithmetic. 
\end{prop}

{\bf Proof.} The non trivial part of the statement is the {\it if} part. Assume $K:=\{0, c, b, a\}$ (where $0 < c < b < a$) is b-arithmetic, then from Proposition \ref{thm:modulus}, the mask polynomial $A(K)=x^a + x^b + x^c + 1$ must have a root $v$ on the unit circle. We want to prove that $v$ is a root of unity. Clearly $v$ satisfies the following equations:
\begin{eqnarray*}
\left( v^{a-b} + 1 \right) v^b &=& - \left(v^c + 1 \right)\\
\left( \overline{v}^{a-b} + 1 \right) \overline{v}^b &=& - \left(\overline{v}^c + 1 \right)
\end{eqnarray*}
where $\overline{v}$ is the complex conjugate of $v$. By multiplying the corresponding sides we get:
$$
v^{a-b} + \overline{v}^{a-b} = v^c + \overline{v}^c
$$

in particular, $v^{a-b}$ and $v^c$ have the same real part. Since they are both on the unit circle, we only have the following two options:
$$
v^{a-b} = v^c \text{ or } v^{a-b}=\overline{v}^c.
$$
that is to say
$$
v^{a-b - c} = 1\text{ or } v^{a-b + c}=1.
$$
If $a-b-c \neq 0$ (from the hypothesis on $a, b,$ and $c$ we cannot have $a - b + c = 0$), we get that $v$ is a root of unity. On the other hand, if $a - b - c=0$, the mask polynomial can be written:
$$
A(K)=x^{b+c} + x^b + x^c + 1=(x^b + 1)(x^c + 1)
$$
and its zeros are all roots of unity. $\square$

\bigskip

{\bf Example 7.} The set $K =\{ 0, 1, 3, 5, 6\}$ is b-arithmetic in $\Z$ but not p-arithmetic. Indeed, again from Proposition \ref{thm:bounded}, we need to prove that the mask polynomial
$$
A(K)=x^6 + x^5 + x^3 + x + 1
$$
has zeros on the unit circle but none of them is a root of unity.

First of all, notice that $A(K)$ is self-reciprocal, therefore, if $\alpha$ is a zero, then $1/\alpha$ is a zero. Moreover, from the resultant, we get that all the zeros are different. Consequently, we can factor $A(K)$ as:
$$
A(K)=(x^2 + a x + 1)(x^2 + b x + 1)(x^2 + c x + 1)
$$
where the solutions to each quadratic polynomial are of the form $\alpha, 1/\alpha$ (notice that 1 and $-1$ are not zeros for $A(K)$) and the coefficients $a,b,$ and $c$ are all distinct. Now, by equating the coefficients we get the following symmetric equations:
\begin{eqnarray*}
a+b+c &=& 1\\
ab + bc + ac &=& -3\\
abc &=& -1.
\end{eqnarray*}

By solving for one variable, say $a$, we get

$$
a^3 - a^2 -3a + 1 = 0,
$$

and by using numerical methods, we get three distinct real solutions, so we can write: $a \approx 2.170$, $b \approx 0.311$, and $c \approx -1.481$. In particular $A(K)$ has 2 real solutions and 4 complex solutions, all of them on the unit circle. Moreover, none of them is a root of unity because $A(K)$ is not divisible by any cyclotomic polynomial. 

\bigskip

 We would like to close this section by pointing out some relevant properties of arithmetic sets in regard to Coven-Meyerowitz Conjecture.
 
 \medskip

(A) $T_1$ implies p-arithmetic: $T_1$ implies in particular that the set $S_K$ is not empty, therefore $P_K(x)$ is divisible by a cyclotomic polynomial; so, from Theorem \ref{thm:equivalent} we have that $K$ is $p$-arithmetic.

 (B) p-arithmetic does not imply $T_1$: for example, consider $K=\{0,3,5,7,9\}$. The mask polynomial equals $1+x^3+x^5+x^7+x^9$, and it can be factored as $$(x^2 - x + 1)(x^7 + x^6 + x^5 + x + 1).$$ The first factor is cyclotomic, but the second is not and $\displaystyle \prod_{p^{\alpha} \in S_K} \phi_{p^{\alpha}}(1) = 1$ while $P_K(1)=5$. 
 
 (C) b-arithmetic does not imply $T_1$: see Example 7.

 (D) p-arithmetic does not imply $T_2$:  for example, consider the set $K=\{0,2,3,5,6, 8\}$. The mask polynomial is  $$P_K(x)=1 + x^2 + x^3 + x^5 + x^6 +x^8=\phi_{4} \cdot \phi_{9}.$$ We have that $\phi_{4}$ and $\phi_{9}$ divide $P_K(x)$, but $\phi_{36}$ does not. It is important to notice that $K$ is not a tile.

 (E) b-arithmetic and not p-arithmetic implies $T_2$ vacuously.

\bigskip

\section{Arithmetic Sets in Free Groups} 
 
 In this section we study the tiles and arithmetic sets in free groups. We will consider only the right tiles, and all the Cayley metrics will be assumed to be left invariant. 
 
 \medskip
 
 \begin{prop}\label{prop:prop1} Every finite connected set of $\mathbb {F}_n, \ n\geq 1$ (more precisely, a finite connected set of vertices in the Cayley graph of $\mathbb {F}_n, \ n\geq 1$ with respect to the standard generating set) is a tile.
\end{prop}

    \medskip

     {\bf Proof.} For $n = 1$ this is obvious. Assume $n\geq 2$, and $\mathbb {F}_n$ is generated by $a_1, \ldots , a_n.$

 \medskip
 
     Let $K$ be a connected set, $\partial K = \{x\in \mathbb {F}_k : d(x,K) = 1\}, \ z_1 \in \partial K, \ |z_1| = min \{|z| : z\notin K \}$.  We will assume that $1\in K$. Clearly, there exists only one $z\in K$ such that $z_1^{-1}z\in \{ a_1, a_1^{-1}, \ldots a_n, a_n^{-1} \}$.  Without loss of generality we may assume that $z_1^{-1}z = a_1$.

\medskip

     Since $K$ is finite, there exists $u_1\in K$ such that $u_1a_1\notin K$. Then, by connectedness of $K$, we have $z_1u_1^{-1}K \cap K = \emptyset $ and $z_1u_1^{-1}K \cup K$ is connected. Then we pick up  $z_2\in \partial (K\cup z_1u_1^{-1}K)$  with  $|z_2| = min \{|z| : z\notin K\cup z_1u_1^{-1}K\}$.  Clearly, $z_2 \notin \partial K \cap \partial z_1u_1^{-1}K$ (i.e. $z_2$ belongs only to one of the sets $\partial K, \partial z_1u_1^{-1}K$) so similarly we may add a third left shift of $K$ disjoint from $K\cup z_1u_1^{-1}K$ such that the three left shifts form a connected subset.  

\medskip

     We continue the process as follows:  suppose the shifts $$g_1K, g_2K, \ldots , g_nK $$ are already chosen such that

  (i) $g_1 = 1$,
  
  (ii) $g_iK\cap g_jK = \emptyset $ for any two different $i, j \in \{1, \ldots ,n\}$,

  (iii) $d(g_{i+1}K, \displaystyle \mathop{\cup }_{j=1}^i g_jK) = 1$  for any $i\in \{1, \ldots n-1\}$.

\medskip

   Then we pick up $z_n\in \partial \displaystyle \mathop{\cup }_{j=1}^n g_jK$ such that $$|z_n| = \mathrm{min} \{|z| : z \notin  \displaystyle \mathop{\cup }_{j=1}^n g_jK\}$$ and since $z_n$ belongs only to one of the sets  $\partial (g_1K) , \ldots , \partial (g_nK)$ \ we may add a new left shift to continue the process. Clearly, the sets $g_1K, g_2K, \ldots $ tile the group.  $\square $

    \medskip

    The converse of Proposition \ref{prop:prop1} also holds, with a slight modification:

    \bigskip

   \begin{prop} \label{prop:prop2} Let $G$ be a finitely generated group with a fixed generating set  $S$. If any connected set (with respect to $S$) of $G$  tiles the group, then $G$ is isomorphic to the free product of some copies of \ $\mathbb{Z}$ \ and \ $\mathbb{Z}/2\mathbb{Z}$. [In particular,  $G$ is virtually free]
\end{prop}

    \medskip

   Proof: We will prove the claim by induction on $|S|$. For $|S| = 1$, \ $G$ is necessarily cyclic, and since a connected set $\{ 0, 1, \ldots n-2\}$ is not a tile of $\mathbb{Z}/n\mathbb{Z}$ for $n\geq 3$,   the claim holds.

\medskip

    Assume $S = \{a_1, a_2, \ldots a_n\}, \ n\geq 2$. For any $i,j\in \{1, \ldots , n\}, i\neq j$, let $G_{ij}$
  be a subgroup of  $G$ generated by $a_i$ and $a_j$. Let also $(V,E)$ be a graph where $V = S, E' = \{ (i, j) : i\neq j$, there exists a path  $r_{ij}(a_i, a_j)$ in  $G_{ij}$ which connects  $a_i$ \ to $a_j$,  and does not pass through the identity element  $1\in G \}$,   $E'' = \{ (i, j) : i\neq j$,  there exists a path  $s_{ij}(a_i, a_j)$   in  $G_{ij}$  which connects  $a_i$ to $a_j^{-1}$, and does not pass through the identity element  $1\in G \}$, $E = E'\cup E''$.

\medskip

    Assume $(V,E)$ is connected. By cyclically permuting  $r_{ij}$ and $s_{ij}$  we obtain paths $t_{ij}$ connecting  $a_i^{-1}$ to $a_j^{-1}$, and  $u_{ij}$ connecting $a_i^{-1}$ to $a_j$ such that none of the $t_{ij}, u_{ij}$  passes through $1\in G$.

\medskip

    Now, let  $F = \displaystyle \mathop{\bigcup }_{1\leq i < j \leq n}R_{ij}\cup S_{ij}\cup T_{ij}\cup U_{ij} $  where 
    $R_{ij}, S_{ij}, T_{ij}, U_{ij}$ are the set of vertices of the paths    $r_{ij}, s_{ij}, t_{ij}, u_{ij}$ respectively (some of these sets could be empty). By construction, $B_1\backslash \{1\} \subset F, \ 1\notin F$, moreover, any connected component of $F$ contains at least two elements. Therefore  $F$ is not a tile.

\medskip

    So $(V,E)$  is not connected: let $(V_1, E_1)$ and $(V_2, E_2)$ be two components with $V = V_1\cup V_2$.  In this case,  $G$ is a free product of two nontrivial subgroups  $G_1, G_2$  where $G_i$ is generated by the elements of $V_i \ (i=1, 2)$,  and we may proceed by induction. $\square $

\bigskip
    
    \begin{prop} \label{prop:last} For $r\geq 2$ and $k\geq 2$, the set $B_{r}\backslash \{1\}$ in the free group $\mathbb{F}_k$ is b-arithmetic but not totally p-arithmetic.
    
    \end{prop} 
    
    \medskip
    
    {\bf Proof.} Let $|.|$ denotes the left invariant Cayley metric with respect to the standard generating set. Let also $K = B_{r}\backslash \{1\}$ and $g_0, g_1, \ldots $ be all elements  of $\mathbb{F}_k$ where the enumeration satisfies the following condition: if $i < j$ then $|g_i|\leq |g_j|$. (i.e. for every $m\geq 0$, we enumerate the elements of the sphere of radius $m$, before starting to enumerate the elements of the sphere of radius $m+1$). 

\medskip

   We will consider the system $A(K)$ and define its solution $(x_g)_{g\in G}$ inductively such that the solution is bounded but not arithmetic. We define $x_g$ for $g\in B_r$ arbitrarily such that $\displaystyle \sum _{g\in K}x_g = 0$ and $x_g\neq 0$ for some $g\in K$. Let $n\geq 1$, and suppose that $x_g$ is defined for all $g\in \displaystyle \mathop{\cup } _{i<n} g_iK$. 

\medskip
   
   Notice that for every $n\geq 1$ the set $g_nK\backslash \displaystyle \mathop{\cup} _{i<n}g_iK$ contains at least two elements, moreover, there exists a positive integer $s = s(n)$ such that  $(\displaystyle \mathop{\cup} _{i<n}g_iK)\cap (B_s\backslash B_{s-1}) = \emptyset $ while $g_nK\cap (B_s\backslash B_{s-1})$ contains at least two elements.
   
\medskip

    Then, for every $n\geq 1$, let $a_n, b_n$ be two distinct elements of the subset $(g_nK\backslash \displaystyle \mathop{\cup}_{i<n}g_iK)\cap (B_s\backslash B_{s-1})$. Then $a_n = g_i, b_n = g_j$ for some $i, j$ (without loss of generality, we may assume that $i < j$), moreover, we define $x_g \in \mathbb{C}$ for all $g\in g_nK$ such that 

\medskip
    
    (i) $x_{a_n}\neq x_k$ for all $k < i$;
    
    (ii) $\sum _{g\in g_nK} x_{g} = 0$.   

\medskip

 Hence, the solution $(x_{g_i})_{i\geq 0}$ is bounded but not periodic. $\square $
 
 \medskip
 
 \begin{rem} Notice that we prove more than the claim of the proposition, namely, we produce a solution $(x_{g_i})_{i\geq 0}$ which is bounded but not periodic. Moreover, the proof of Proposition \ref{prop:last} holds also for an arbitrary finite subset $K$ where (i) $1\notin K$; (ii) $B_2(1)\backslash \{1\}\subseteq K$; and (iii) $\{1\}\cup K$ is connected. Notice that any such set $K$, necessarily,  is not a tile. On the other hand, sets like this can easily be p-arithmetic. To produce such an example, let $k = 2$, and $K$ be a subset of $B_3(1)$ satisfying conditions (i)-(iii). Let also $|K| = 24$ (so $|K\backslash B_2(1)| = 8$). Define a vector ${\bf x} = (x_g)_{g\in \mathbb{F}_2}$ by letting $x_g = 1$ whenever $|g|$ is odd, and $x_g = -1$ whenever $|g|$ is even. Then ${\bf x}$ is a periodic non-degenerate solution of $A(K)$.     
 \end{rem}

\end{document}